\newcommand{\bull}{\vrule height .9ex width .8ex depth -.1ex}
 \newcommand{\ppp}{\hfill $\bull$ }
 \documentclass[11pt]{article}
 \author{ Mohammed-Larbi Labbi}
   \title{Variational Properties of the Gauss-Bonnet Curvatures}
   \date{}
   \usepackage{amsmath}
   \usepackage{amscd}
\newtheorem{theorem}{Theorem}[section]
\newtheorem{corollary}[theorem]{Corollary}
\newtheorem{lemma}[theorem]{Lemma}

\newtheorem{proposition}[theorem]{Proposition}
\newtheorem{definition}[theorem]{Definition}
\begin{document}
   \maketitle
   \begin{abstract} The Gauss-Bonnet curvature of order $2k$ is a generalization to higher
dimensions of the Gauss-Bonnet integrand in dimension $2k$,
as the scalar curvature generalizes the
two dimensional Gauss-Bonnet integrand.\\
In this paper, we evaluate the first variation of the integrals of these curvatures seen as
 functionals on the space of all Riemannian metrics on the manifold under consideration. 
An important property of this derivative is that it depends only on the curvature tensor 
and not 
on its covariant derivatives.
We show that the critical
points of this functional once restricted to metrics with unit  volume  
are generalized Einstein metrics
and once restricted to a pointwise conformal class of metrics are metrics with constant
Gauss-Bonnet curvature.
\end{abstract}
   \par
 \section{Introduction and Statement of the Results}
 Let $M$ be a compact smooth (oriented) manifold of dimension $n$, and let
  ${\cal M}$ be the space of smooth
 Riemannian metrics on $M$  endowed with a natural  $L^2$-Sobolev norm.
 This allows us to speak about differentiable functionals
 ${\cal M}\rightarrow {\bf R}$.
A functional $F:{\cal M}\rightarrow {\bf R}$ is called Riemannian if it
 is invariant under the action of the diffeomorphism group. We say that $F$
 has a gradient at $g$ if there exists a symmetric tensor $a\in {\cal C}^1$
  such that
 for every symmetric tensor $h\in {\cal C}^1$ we have
 $$F_g'h={d\over dt}\mid_{t=0}F(g+th)=<a,h>,$$
where ${\cal C}^1$ denotes the space of symmetric tensors in
$\Lambda^{*}M \otimes\Lambda^{*}M$ and
$<,>$ is the integral inner product.
\par\medskip\noindent
 A classical  Riemannian functional is  the total scalar curvature
 $$S(g)=\int_M{\rm scal}\mu_g,$$
 where ${\rm scal}$ denotes the scalar curvature function of $g$
  and $\mu_g$ is the
 volume element of $g$.  An important point about this functional $S$,
 is that its critical points, when restricted to
  ${\cal M}_1=\{ g\in {\cal M}:{\rm vol}(g)=1\}$,
 are Einstein metrics. Also, its gradient is the Einstein tensor,
 precisely
 $$S_g'h=<{1\over 2}{\rm scal}\, g-{\rm Ric},h>.$$
 Where  ${\rm Ric}$ denotes the Ricci tensor.\\
 A natural generalization of  the functional $S$ is the Riemannian
  functional
  $$H_{2k}(g)=\int_Mh_{2k}\mu_g,$$
  where, for each $1\leq k\leq n/2$, $h_{2k}$ is the Gauss-Bonnet curvature of order $2k$.This
curvature is determined by the complete contraction of the Gauss-Kronecker tensor of order $k$.
 Furthermore, These curvatures coincide with the intrinsic curvature
   invariants of $(M,g)$ which appear in
the well known tube formula of H. Weyl \cite{Wey}, see section 2 below for precise definitions.
 \par
 For $k=1$, $H_2=S/2$ is one  half  the total scalar curvature functional.
  Also, if the dimension $n$ of $M$ is even, then $H_n$ does not depend on
  the metric. It is, up to a constant,
  the Euler-Poincar\'e characteristic of $M$.\par
  \medskip\noindent
  Marcel Berger proved in \cite{Berger} that the gradient of $H_4$,
  like the gradient of $S$, depends only on the curvature tensor $R$
   and does not include
  its covariant derivatives. The expression of the gradient  he obtained was  complicated  
and  hardly generalizable to  higher
  $H_{2k}$. So he asked the following two questions:
  \begin{itemize}
  \item
   Does the above phenomena remain true for all  higher $H_{2k}$?.
   \item Characterize the critical Riemannian metrics for the functional $H_{2k}$.
   \end{itemize}
   \par\medskip\noindent
  In this paper, we completely answer  the first
  question, and give partial answers to the second one. The main result of this paper 
is the following:\par\medskip
\noindent
  {\bf Main Theorem.}\,
  {\sl For every compact n-dimensional Riemannian manifold $(M,g)$,
  and for every integer $k$ such that $2\leq 2k\leq n$,
the functional $H_{2k}$ is differentiable, and at $g$ we have
  $$H_{2k}'h={1\over 2}\langle h_{2k}g- {1\over (2k-1)!}c^{2k-1}R^k,h\rangle.$$
Where $R^k$ denotes the exterior product of the Riemann curvature tensor $R$ with
 itself $k$-times in the ring of curvature structures, $c$ is the
 contraction map and $\langle, \rangle$ is the integral scalar product.}\\
\par\noindent\bigskip
Remark that, for $k=1$, the main theorem shows that
$$H_{2}'h={1\over 2}<h_2g-
cR,h>={1\over 2}<{{\rm scal}\over 2} g- {\rm Ric},h>.$$
So that we recover the above formula   about the total scalar curvature. Also,
 in the case $2k=n$, we have (see (\ref{justif}) below):
$$H_{n}'h={1\over 2}<h_{n}g-
{1\over (n-1)!}c^{n-1}R^k,h>=0.$$
This is not a surprise, because $H_n$ does not depend on the metric
 by the Gauss-Bonnet theorem.
\par\medskip
\noindent
Let us note here that the main theorem was established earlier by
 David Lovelock in a not
well known paper to the mathematicians \cite{lovelock}. His proof is based
on classical tensor analysis.
Our proof of the main theorem is simple and coordinate free. The key point of our proof
is  that it is possible to write the Gauss-Bonnet curvatures   as
 exterior products of the metric $g$  with
  the Riemann
 curvature tensor $R$ (see the definition below for a precise formulation),
  then one can get the
 desired derivative using the power rule of differentiation
  and stokes' theorem.\\
The main theorem shall be proved in section 4. In section  2,
 we recall some useful facts
about the ring of curvature structures from  \cite{Labbi}. In section 3,
 we show that many of the  classical
results of Hodge theory  on differential forms can be naturally extended
 to the context of double forms. We consider here
  only those results which
 shall be used later in this paper.
Then, we define and study  some operators on double forms
which will play
a key role in the proof of the main theorem.
\par\noindent In section 5, we study the critical metrics
of the
functional $H_{2k}$, when restricted to a normalized conformal class of some
metric, and we prove that  they are  metrics with
constant ($2k$)-Gauss-Bonnet curvature. It is then natural to ask
whether in each conformal class of a Riemannian metric,  on a
smooth compact manifold of dimension $n>2k$, there exists a metric
with constant ($2k$)-Gauss-Bonnet curvature. That is  a natural
generalization of the
famous Yamabe problem.\par
\noindent
Finally, in section 6 we examine some properties of  the critical metrics of the functional
 $H_{2k}$ in
the space of all Riemannian metrics with volume 1. These are generalized Einstein metrics.
 
 \section{Preliminaries}
Let $(M,g)$ be a compact smooth (oriented) Riemannian manifold of dimension n. We denote by
 $\Lambda^{*}M=\bigoplus_{p\geq 0}\Lambda^{*p}M$ the ring of differential
  forms
 on $M$.  Considering the tensor product over the ring of  smooth functions,
  we define
 ${\cal D}= \Lambda^{*}M\otimes \Lambda^{*}M=\bigoplus_{p,q\geq 0}
  {\cal D}^{p,q}$ where $  {\cal D}^{p,q}= \Lambda^{*p}M \otimes\Lambda^{*q}M$.
  It is a graded associative  ring and it is  called the ring of
   double forms on $M$. The exterior product in ${\cal D}$, sometimes called the 
Kulkarni-Nomizu product,  will be denoted
   by a dot., this shall be omitted whenever possible. \par\noindent
   The ring of curvature structures on $M$ is the ring ${\cal C}=\sum_{p\geq 0}
   {\cal C}^p$ where ${\cal C}^p$ denotes symmetric elements in
   ${\cal D}^{p,p}$.\par\medskip
The standard inner product $<,>$ on $\Lambda^{*p}M$ and the Hodge star
operator $*$ extend
in a natural way    to ${\cal D}$ (we assume here that the manifold is orientable). These were used in
\cite{Labbi} to study several properties of this ring. In particular,
 the following   relations are proved for all
 $\omega,\omega_1,\omega_2\in {\cal D}$:
 \begin{equation}\label{g=starc}
g\omega=*c*\omega,
\end{equation}
\begin{equation}\label{conjugg}
 <g\omega_1,\omega_2>=<\omega_1,c\omega_2>.
 \end{equation}
 Where $c$ denotes the contraction map.\\
 Also, for all  $\omega,\theta\in D^{p,q}$, we have
\begin{equation}\label{IPHodge}
<\omega,\theta>=*(\omega.*\theta)=(-1)^{(p+q)(n-p-q)}*(*\omega.\theta),
\end{equation}
\begin{equation}\label{starsquare}
**\omega =(-1)^{(p+q)(n-p-q)}\omega.
 \end{equation}
{\sc Remark.} A double form $\omega\in {\cal D}^{p,q}$ can be seen as a symmetric bilinear
 form acting on $p$-vectors. Under
this identification one can check easily that 
\begin{equation*}
*\omega(.,.)=(-1)^{(p+q)(n-p-q)}\omega(*.,*.).
\end{equation*}
The minus sign appears in fact  because for a usual $p$-form $\theta$ 
and for the usual Hodge star operator we have $(*\theta)(.)=(-1)^{p(n-p)}\theta(*.)$.\\

In the following, we shall denote by $\omega^q$,  the product of $\omega$ with itself
q-times in the ring ${\cal C}$.\par\medskip\noindent
Let us recall the following definitions:
\begin{definition}{\bf \cite{Labbi}} The $(p,q)$-curvature, denoted $s_{(p,q)}$,
  for $1\leq q\leq {n\over 2}$ and $0\leq p\leq n-2q$, is the sectional
  curvature
  of the following $(p,q)$-curvature tensor
  \begin{equation}\label{spq:def}
  R_{(p,q)}={1\over (n-2q-p)!}*\bigl( g^{n-2q-p}R^q\bigr)
  \end{equation}
 In other words, for a tangent $p$-plane $P$, $s_{(p,q)}(P)$ is the sectional curvature of the tensor
  ${1\over (n-2q-p)!}g^{n-2q-p}R^q$ at the orthogonal complement
of $P$.
\end{definition}
Note that the tensors $R_{(p,q)}$ satisfy the first Bianchi identity and they are
divergence free.\par\medskip\noindent
Here in this paper, we are mainly interested in the following  special
cases, the $(0,q)$ and $(1,q)$-curvatures:
\begin{definition}{\bf ( \cite{Labbi}, \cite{Labbi2})} Let $q$ be a positive integer such that $2\leq 2q\leq n$.
\begin{enumerate}
\item The $(2q)$-Gauss-Bonnet curvature,   denoted $h_{2q}$, is
the $(0,q)$-curvature. That is the function defined on $M$ by
\begin{equation}
h_{2q}=s_{(0,q)}={1\over (n-2q)!}*\bigl( g^{n-2q}R^q\bigr).
\end{equation}
\item The $(2q)$-Einstein-Lovelock tensor, denoted $T_{2q}$, is defined to be the $(1,q)$-curvature
tensor, that is
\begin{equation}\label{t2q:def}
T_{2q}=*{1\over (n-2q-1)!}g^{n-2q-1}R^q.
\end{equation}
\end{enumerate}
If $2q=n$, we set $T_n=o$.
\end{definition}
Recall the following properties of these curvatures  \cite{Labbi}
\begin{equation}\label{h2q:def}
h_{2q}= {1\over (n-2q)!}*\bigl( g^{n-2q}R^q\bigr)={1\over (2q)!}c^{2q}R^q,
\end{equation}
that is the complete contraction of $R^q$. Also,
\begin{equation}\label{t2q:property}
T_{2q}={1\over (2q)!}c^{2q}R^qg-{1\over (2q-1)!}c^{2q-1}R^q=h_{2q}g-
{1\over (2q-1)!}c^{2q-1}R^q.
\end{equation}
For $q=1$ we recover the usual Einstein tensor $T_2={1\over 2}c^2Rg-cR$.
 For $2q=n$, we have
$R^q=\lambda {g^n\over n!} $ for some constant $\lambda$, hence
\begin{equation}\label{justif}
h_{2q}g- {c^{2q-1}R^q\over (2q-1)!}={1\over n!}
c^n(\lambda {g^n\over n!})g-
{1\over (n-1)!}c^{n-1}(\lambda {g^n\over n!})=\lambda g-\lambda g=0.
\end{equation}
This justifies our definition for $T_n$.\par\smallskip\noindent
Note that $c^{2q-1}R^q$ can be considered as a generalization of
the Ricci
 curvature tensor.\\
{\sc Remark.}
Let us remark here some analogies between the Einstein-Lovelock tensors and the
usual Einstein tensor. Recall that, the former is the main linear combination
of the metric tensor $g$ and its Ricci curvature to be divergence free. The same property is true 
for  Einstein-Lovelock tensors if we
substitute the  generalized Ricci curvatures $c^{2q-1}R^q$  to the usual one.\\
Another similarity between these tensors can be noticed at the level of their sectional
curvatures: The full contraction of the Riemann curvature tensor $R$ in the directions
orthogonal to a given direction $v$ produces the sectional curvature of the
usual Einstein tensor, that is $T_2(v,v)$. Doing the same operation but for the Gauss-Kronecker tensors
$R^q$ one generates the sectional curvatures of the Einstein-Lovelock tensors, that is 
$T_{2q}(v,v)$.\\
Finally, let us recall the following property \cite{Labbi} which provides another analogy:
\begin{equation}\label{tracet2q}
  {\rm trace}\, T_{2k}=(n-2q)h_{2q}.\end{equation}

\section{The Second Bianchi Map and other Differential Operators }
The second Bianchi sum, denoted  $D$, maps ${\cal D}^{p,q}$ into  ${\cal D}^{p+1,q}$.
 For $\omega \in  {\cal D}^{p,q}$, we have
$$(D\omega)(x_1\wedge...\wedge x_{p+1},y_1\wedge...\wedge y_{q})=
\sum_{j=1}^{p+1}(-1)^j{\nabla_{x_j}\omega}(x_1\wedge...\wedge
\hat{x_j}\wedge ... x_{p+1}, y_1\wedge...\wedge y_{q}),$$
where $\nabla$ denotes the covariant differentiation with respect to
 the metric $g$. \par
 If we identify  double forms with vector valued differential forms, then  $D$ coincides
 with the operator of exterior differentiation of vector valued differential forms \cite{Besse}.
 In particular, 
the restriction of $D$ to ${\cal D}^{p,0}$ coincides with $-d$, where $d$ is the
 usual  exterior derivative on $p$-forms.
 A second possible extension of $d$ is the adjoint second Bianchi sum:
 $$\tilde{D}:{\cal D}^{p,q} \rightarrow {\cal D}^{p,q+1}$$
 defined for $\omega\in{\cal D}^{p,q}$ by
$$(\tilde{D}\omega)(x_1\wedge...\wedge x_{p},y_1\wedge...\wedge y_{q+1})=
\sum_{j=1}^{q+1}(-1)^j{\nabla_{y_j}\omega}(x_1\wedge...\wedge x_{p},
 y_1\wedge... \wedge \hat{y_j}\wedge ...\wedge y_{q+1}).$$
 Note that in general we have neither $D^2=0$ nor ${\tilde D}^2=0$.
 The composition of these two operators is the operator

 $$D\tilde D:{\cal D}^{p,q} \rightarrow {\cal D}^{p+1,q+1}$$
\begin{equation}\label{dtilde:d}
\begin{split}
(&D\tilde{D}\omega)(x_1\wedge...\wedge x_{p+1},y_1\wedge...\wedge y_{q+1})=\\
&\sum_{i=1}^{p+1}\sum_{j=1}^{q+1}(-1)^{i+j}(\nabla^2_{x_iy_j}\omega)
(x_1\wedge...\wedge \hat{x_i}\wedge ...\wedge x_{p+1},
 y_1\wedge... \wedge \hat{y_j}\wedge ...\wedge y_{q+1}).
\end{split}
\end{equation}
 Remark that the restriction of $D\tilde D$ to $D^{0,0}$ is the usual
  Hessian
 operator on functions, precisely
 $$D\tilde D(f)(x,y)=\nabla^2_{x,y}f.$$
 Also, note that for $h\in D^{1,1}$, we have
 \begin{equation}
 \label{DDh:form}D\tilde D h(x\wedge y,z\wedge u)=\nabla^2_{xz}h(y,u)-\nabla^2_{xu}h(y,z)-\nabla^2_{yz}h(x,u)
 +\nabla^2_{yu}h(x,z).
 \end{equation}
 Similarly, one can  also consider the differential operator
  $$\tilde D D:{\cal D}^{p,q} \rightarrow {\cal D}^{p+1,q+1}$$
\begin{equation*}
\begin{split}
(&\tilde{D} D\omega)(x_1\wedge...\wedge x_{p+1},y_1\wedge...\wedge
y_{q+1})=\\
&\sum_{i=1}^{p+1}\sum_{j=1}^{q+1}(-1)^{i+j}(\nabla^2_{y_jx_i}\omega)
(x_1\wedge...\wedge \hat{x_i}\wedge ...\wedge x_{p+1},
 y_1\wedge... \wedge \hat{y_j}\wedge ...\wedge y_{q+1}).
\end{split}
\end{equation*}
If $\omega$ is a symmetric double form, $D\tilde D\omega$ and
$\tilde DD\omega$ are not necessarily symmetric. Nevertheless, it is true that
$$(\tilde{D} D\omega)(x_1\wedge...\wedge x_{p+1},y_1\wedge...\wedge y_{q+1})
=(D\tilde{D} \omega)(y_1\wedge...\wedge y_{q+1},x_1\wedge...\wedge x_{p+1}).$$
Consequently, the following operator is well defined and will play an
 important role
in the proof of the main theorem
$$D\tilde D+\tilde D D:{\cal C}^p\rightarrow {\cal C}^{p+1}$$
This operator can be considered as a natural generalization, to the ring
${\cal C}^p$,  of the usual
Hessian operator on functions.
\par\medskip\noindent
On the other hand, it is easy to check that for $\omega\in D^{p,q}$ and $\theta\in D^{r,s}$
 we have
 \begin{equation}\label{dom:thet}
 \begin{split}
 D(\omega.\theta)&=D\omega.\theta +(-1)^p\omega.D\theta,\\
 \tilde D(\omega.\theta)&=\tilde D\omega.\theta +(-1)^q\omega.\tilde D\theta.
 \end{split}
 \end{equation}
 The operator $\delta=c\tilde D+\tilde D c$ was defined by Kulkarni
 \cite{Kulk} as a natural
 generalization of the classical $\delta$ operator.
 Using the Hodge star operator we shall extend some classical results of
  Hodge theory to  double forms as follows:
  \begin{proposition} If $ *$  denotes the generalized Hodge star operator on ${\cal D}$, then
  we have
   \begin{equation}\label{stard:star}\delta \omega  =(-1)^{(p+q)(n-p-q)}*D* \omega \end{equation}
for every $(p,q)$-form $\omega$ such that $p\geq 1$.
\end{proposition}
{\bf Proof.}  Let $ \omega\in D^{p,q}$, and $(x_i)$, $(y_j)$ orthonormal
  vector fields about $m\in M$ such that
   $(\nabla_{x_i})_m=(\nabla_{y_j})_m=0$. Then at $m$ we have
  \begin{equation*}
  \begin{split}
  (& D*\omega)(x_1\wedge ...\wedge x_{n-p+1},y_1\wedge ...\wedge y_{n-q})\\
  =&\sum_{j=1}^{n-p+1}(-1)^j \nabla_{x_j}(*\omega)
  (x_1\wedge ...\wedge {\hat x}_j \wedge ...\wedge x_{n-p+1},y_1\wedge ...
  \wedge y_{n-q})\\
  =&\sum_{j=1}^{n-p+1}(-1)^j \nabla_{x_j}(*\omega
  (x_1\wedge ...\wedge {\hat x}_j \wedge ...\wedge x_{n-p+1},y_1\wedge ...
  \wedge y_{n-q}))\\
  =&\sum_{j=1}^{n-p+1}(-1)^{j+j-1+(p+q)(n-p-q)} \nabla_{x_j}\left(\omega
  (x_j\wedge *(x_1\wedge  ...\wedge x_{n-p+1}),*(y_1\wedge ...
  \wedge y_{n-q})\right)\\
=&-\sum_{j=1}^{n-p+1}(-1)^{(p+q)(n-p-q)}\nabla_{x_j}(\omega)
  (x_j\wedge * (x_1\wedge ...\wedge x_{n-p+1}),*(y_1\wedge ...
  \wedge y_{n-q}))\\
  =&(-1)^{(p+q)(n-p-q)}\delta \omega ( * (x_1\wedge ...\wedge x_{n-p+1}),*(y_1\wedge ...
  \wedge y_{n-q}))\\
  =& (-1)^{(p+q)(n-p-q)}(-1)^{(p+q-1)(n-p-q+1)}*\delta \omega(x_1\wedge ...
\wedge x_{n-p+1},y_1\wedge ...
  \wedge y_{n-q}).
\end{split}
  \end{equation*}
Therefore $D*\omega=(-1)^{n+1}*\delta \omega$.
  \ppp
\begin{proposition} With respect to the integral scalar product, the operator
\begin{equation}\label{adjoint:D}
(-1)^{n+p}\delta:D^{p+1,q}\rightarrow D^{p,q}
\end{equation}
 is the formal
adjoint of the operator D.\end{proposition}
{\bf Proof.} Let $\omega_1 \in D^{p,q}$ and $\omega_2\in D^{p+1,q}$ then
\begin{equation*}
\begin{split}
D(\omega_1.*\omega_2)&=D\omega_1.*\omega_2+(-1)^p\omega_1.D*\omega_2\\
=& D\omega_1.*\omega_2+(-1)^p (-1)^{n+1}\omega_1.*\delta\omega_2\\
=& **\{ (D\omega_1).*\omega_2\}+(-1)^{n+p+1}**\{\omega_1.*\delta\omega_2\}\\
=&*\biggl( \langle D\omega_1,\omega_2\rangle +(-1)^{n+p+1}
\langle \omega_1,\delta\omega_2 \rangle \biggr).\\
\end{split}
\end{equation*}
Applying the generalized Hodge operator to both sides of the previous equation we get
\begin{equation*}
-\delta(*(\omega_1.*\omega_2))=<D\omega_1,\omega_2>+(-1)^{n+p+1}
<\omega_1,\delta\omega_2>.
\end{equation*}
Note that, the integral of the left hand side is zero by stokes' theorem. This completes
the proof.\ppp
\par\medskip\noindent
In the same way, one can prove without difficulties that if $\tilde \delta=cD+Dc$ then for every
$(p,q)$-form $\omega$ with $q\geq 1$ we have
\begin{equation}\label{stardtilde:star}
\tilde \delta \omega=(-1)^{(p+q)(n-p-q)}*\tilde D *\omega,\end{equation}
and also that
\begin{equation}\label{adjoint:Dtilde}
(-1)^{n+q}\tilde\delta:D^{p,q+1}\rightarrow D^{p,q}
\end{equation}
is the formal adjoint of $\tilde D$ with respect to the integral scalar product.
\begin{corollary}\begin{enumerate}
\item  For $p,q\geq 1$, The operators $\tilde\delta \delta$ and
 $\delta\tilde\delta $ send
$D^{p,q}$ to $D^{p-1,q-1}$ and respectively satisfy
\begin{equation}\label{adjoint:tildedd}
\tilde\delta \delta=(-1)^{(p+q)(n-p-q)}*\tilde D D*,\quad \quad 
\delta \tilde\delta =(-1)^{(p+q)(n-p-q)}*D \tilde D *.
\end{equation}
Furthermore,  with respect to the integral scalar product, they are respectively 
the formal adjoints 
of the operators $(-1)^{p+q}D\tilde D$ and $(-1)^{p+q}\tilde D D$.
\item For $p,q\geq 1$, The operator $\tilde\delta \delta +\delta\tilde\delta $ sends 
$D^{p,q}$ to $D^{p-1,q-1}$ and  satisfies
\begin{equation}\label{adjoint:dtilded}
\tilde\delta \delta +\delta\tilde\delta=(-1)^{(p+q)(n-p-q)}*(\tilde D D+D \tilde D) *.
\end{equation}
Furthermore,  with respect to the integral scalar product, it is  
the formal adjoint
of the operator $(-1)^{p+q}(D\tilde D+\tilde D D)$.
\end{enumerate}
\end{corollary}
{\bf Proof.} It is a direct consequence of (\ref{adjoint:D}) and
(\ref{adjoint:Dtilde}). \ppp
\par\medskip\noindent
Remark that since $Dg=0$ and $DR=0$, then $D(g^pR^q)=0$. However the tensors $R_{(p,q)}$
 do not in general satisfy the second Bianchi identity. Nevertheless, they 
are divergence free, In fact,
$$\delta(*g^pR^q)=*D(g^pR^q)=0.$$
This fact is used to prove the following proposition:
\begin{proposition}[Schur's Theorem]
 Let $p,q\geq 1$.
 If at every point $m\in M$ the $(p,q)$-curvature is a constant
 (that is on the fiber at m),
 then it is a constant.\end{proposition}
 {\bf Proof.} Recall that at each  $m\in M$ we have
 $$s_{(p,q)}=\lambda  \Leftrightarrow R_{(p,q)}=
 \lambda{g^{p}\over p!},$$
where $\lambda=\lambda(m)$ is  constant at $m$.\par\noindent
 Since the tensors $R_{(p,q)}$ are divergence free, then
 $$\delta \bigl(\lambda{g^{p}\over p!}\bigr)=0,$$
 and therefore, (see \cite{Labbi})
$$D\bigl(*\lambda{g^{p}\over p!}\bigr)=D\bigl(\lambda{g^{n-p}\over
(n-p)!}\bigr)=0.$$
Consequently $d\lambda =0.$
 This completes the proof.\ppp
 \par\medskip\noindent

Another operator which also
 will play a central role in the
proof of the main theorem is defined as follows:
\par\smallskip\noindent
For each $h\in {\cal C}^1$, we define the operator $F_h:{\cal C}^p\rightarrow
{\cal C}^p$ as follows. Let $m\in M$ and $\{e_1,...,e_n\}$ be an orthonormal basis
 of $T_mM$ diagonalizing $h$ , then its value on  basis elements is
\begin{equation}\label{fh:om}
\begin{split}
F_h&\omega(e_{i_1}\wedge ...\wedge e_{i_p},e_{j_1}\wedge ...\wedge e_{j_p})=\\
&\biggl(\sum_{k=1}^ph(e_{i_k},e_{i_k})+\sum_{k=1}^ph(e_{j_k},e_{j_k})\biggr)\omega
(e_{i_1}\wedge ...\wedge e_{i_p},e_{j_1}\wedge ...\wedge e_{j_p}).
\end{split}
\end{equation}
It is not difficult to see that if $\omega$
satisfies the first Bianchi identity then  so does $F_h(\omega)$.
Below we shall  prove some useful  properties of this operator.
\begin{proposition}\label{fh:derivation} For all $\omega\in {\cal C}^p$ and $\theta \in {\cal C}^q$
we have
$$F_h(\omega.\theta)=F_h(\omega).\theta+\omega.F_h(\theta).$$
That is, $F_h$ acts by derivations on ${\cal C}$. In particular we have
\begin{equation}\label{fhomk}
F_h(\omega^k)=k\omega^{k-1}F_h(\omega).
\end{equation}
\end{proposition}
{\bf Proof.}  Let $m\in M$ and $\{e_1,...,e_n\}$ be an orthonormal basis
 of $T_mM$ diagonalizing $h$. Let $\{i_1, ..., i_{p+q}\}$ and
  $\{j_1, ..., j_{p+q}\}$
be arbitrary subsets of $\{1, ...,n\}$ both with $p+q$ elements,
 then at $m$ we have
\begin{equation*}
\begin{split}
&\omega .F_h(\theta)(e_{\scriptstyle i_1}\wedge ...\wedge
e_{\scriptstyle i_{p+q}},
e_{\scriptstyle j_1}\wedge ...\wedge e_{\scriptstyle j_{p+q}})=\\
&\frac{1}{ (p!)^2(q!)^2}
\biggl(\sum_{\scriptstyle \sigma ,\rho \in S_{p+q}}
\epsilon(\sigma)\epsilon(\rho)
\omega(e_{\scriptstyle i_{\sigma(1)}}\wedge...\wedge
e_{\scriptstyle i_{\sigma(p)}},
e_{\scriptstyle j_{\rho(1)}}
\wedge...\wedge e_{\scriptstyle j_{\rho(p)}})\\
&{\phantom {mmmmmmmmmm}}
F_h(\theta)(e_{\scriptstyle i_{\sigma(p+1)}}\wedge...\wedge
 e_{\scriptstyle i_{\sigma(p+q)}},
e_{\scriptstyle j_{\rho(p+1)}}
\wedge...\wedge e_{\scriptstyle j_{\rho(p+q)}})\biggr)\\
&=\frac{1}{(p!)^2(q!)^2}
\biggl( (\sum_{\scriptstyle \sigma ,\rho\in S_{p+q}}
\epsilon(\sigma)\epsilon(\rho)
\omega(e_{\scriptstyle i_{\sigma(1)}}\wedge...\wedge
e_{\scriptstyle i_{\sigma(p)}},
e_{\scriptstyle j_{\rho(1)}}
\wedge...\wedge e_{j_{\rho(p)}})\\
&\phantom{m}
\bigl(\sum_{k=1}^{p+q}h(e_{\scriptstyle i_k},e_{\scriptstyle i_k})+
\sum_{k=1}^{p+q}h(e_{\scriptstyle j_k},
e_{\scriptstyle j_k})
-\sum_{k=1}^{p}h(e_{\scriptstyle i_{\scriptstyle\sigma(k)}},
e_{\scriptstyle i_{\scriptstyle \sigma(k)}})-
\sum_{k=1}^{p}h(e_{\scriptstyle j_{\rho(k)}},
e_{\scriptstyle j_{\rho(k)}})\bigr)\\
&\phantom{mmm}
\theta(e_{\scriptstyle i_{\scriptstyle\sigma(p+1)}}\wedge...\wedge
 e_{\scriptstyle i_{\scriptstyle\sigma(2p)}},
e_{\scriptstyle j_{\scriptstyle\rho(q+1)}}
\wedge...\wedge e_{j_{\scriptstyle\rho(2q)}})\biggr)\\
&=\{F_h(\omega.\theta)-F_h(\omega).\theta\}(e_{\scriptstyle i_1}\wedge ...
\wedge
 e_{\scriptstyle i_{p+q}},
e_{\scriptstyle j_1}\wedge ...\wedge e_{\scriptstyle j_{p+q}}).
\end{split}
\end{equation*}
This completes the proof of the proposition.\ppp

\begin{proposition} The operator $F_h$ is  self adjoint, precisely
for all $\omega,\theta\in {\cal C}$ we have
\begin{equation}\label{fh:autoadjoint}
 <F_h(\omega),\theta>=<\omega,F_h(\theta)>.
 \end{equation}
 \end{proposition}
{\bf Proof.} Let $\{e_1,...,e_n\}$ be an orthonormal basis diagonalizing $h$
 at $m\in M$,
let $\{i_1, ..., i_p\}$ and
  $\{j_1, ..., j_p\}$
be arbitrary subsets of $\{1, ...,n\}$ both with $p$ elements,
 then at $m$  we have
\begin{equation*}
\begin{split}
&F_h(\omega)(e_{\scriptstyle i_1}\wedge ...\wedge e_{\scriptstyle i_{p}},
e_{\scriptstyle j_1}\wedge ...\wedge e_{\scriptstyle j_{p}})
\theta (e_{\scriptstyle i_1}\wedge ...\wedge e_{\scriptstyle i_{p}},
e_{\scriptstyle j_1}\wedge ...\wedge e_{\scriptstyle j_{p}})\\
&=\bigl(\sum_{k=1}^{p}h(e_{i_k},e_{i_k}) + h(e_{j_k},e_{j_k})\bigr)
\omega(e_{\scriptstyle i_1}\wedge ...\wedge e_{\scriptstyle i_{p}},
e_{\scriptstyle j_1}\wedge ...\wedge e_{\scriptstyle j_{p}})\\
&\phantom{mmmmmmmm} \theta (e_{\scriptstyle i_1}\wedge ...\wedge e_{\scriptstyle i_{p}},
e_{\scriptstyle j_1}\wedge ...\wedge e_{\scriptstyle j_{p}})\\
&=\omega(e_{\scriptstyle i_1}\wedge ...\wedge e_{\scriptstyle i_{p}},
e_{\scriptstyle j_1}\wedge ...\wedge e_{\scriptstyle j_{p}})
F_h(\theta) (e_{\scriptstyle i_1}\wedge ...
\wedge e_{\scriptstyle i_{p}},e_{\scriptstyle j_1}\wedge ...\wedge
 e_{ \scriptstyle j_{p}}).
\end{split}
\end{equation*}
The proposition results immediately after taking the corresponding sums.\ppp

\par\medskip\noindent
The following properties can be checked without difficulties:
\par\medskip\noindent
\begin{proposition}\label{extra} The following are true about the operators $F_h$:
\begin{enumerate}
\item If $\omega \in {\cal C}^p$ then $c^p\bigl( F_h(\omega)\bigr)=2p
\langle c^{p-1}\omega, h \rangle.$
\item If $\omega \in {\cal C}^p$ then $F_g(\omega)=2p\omega$.
\item If $p\geq 1$, then $F_h(g^p)=2pg^{p-1}h.$
\item If $h,k\in {\cal C}^1$ and $\bar h,\bar k, \overline{ F_h(k)}$ denote respectively
the associated linear operators on the tangent space, then
\begin{equation*}
\overline{F_h(k)}=\bar h o \bar k+ \bar k o \bar h.
\end{equation*} 
\item If $\omega\in {\cal C}^2$ and $ h\in {\cal C}^1$, then for all
 tangent vectors $x,y,z,u$ we
have:
$$F_h(\omega)(x\wedge y,z\wedge u)=h(\omega(x,y)z,u)-h(\omega(x,y)u,z)+h(\omega(z,u)x,y)-
h(\omega(z,u)y,x),$$
where in the right hand side, $\omega$ was considered  as a $(3,1)$-tensor by the mean of the 
metric $g$.
\item If $\omega\in{\cal C}^n$, then  $F_h(\omega)=2({\rm tr}_gh)\omega$. \end{enumerate}
\end{proposition}
{\sc Remark.} The operator $F_h$ can also alternatively be defined by declaring that it acts by derivations
on ${\cal C}$ and that its restriction to ${\cal C}^1$ is the symmetric multiplication
by $h$ as in the property 4 of  the previous proposition.
\section{Proof of the Main Theorem}
The proof is based on the following two lemmas. The first lemma asserts that
the directional derivative of the Riemann curvature tensor when considered
as a symetric double form is the sum of an "exact" double form and a term which
depends linearly on the curvature:\par\medskip\noindent
\begin{lemma} The derivative of the Riemann curvature structure
$R\in {\cal C}^2$ in the direction of $h\in {\cal C}^1$ is given by
\begin{equation}\label{derivative:R}
R'h={-1\over 4}(D\tilde D+\tilde D D)(h)+{1\over 4}F_h(R).
\end{equation}
\end{lemma}
{\bf Proof.}
Recall that, the directional derivatives $\nabla'h$, $R'h$ at $g$
 of the  Levi-Civita
connection and the  $(3,1)$-Riemann curvature tensor
are respectively given by  \cite{Besse}
$$g(\nabla'h(x,y),z)={1\over 2}\{\nabla_xh(y,z)+\nabla_yh(x,z)-\nabla_zh(x,y)\},$$
$$R'h(x,y)z=(\nabla_y\nabla'h)(x,z)-(\nabla_x\nabla'h)(y,z).$$
Therefore the derivative of $R$, seen as a double form
$R\in \Lambda^{*2}M\otimes \Lambda^{*2}M$, in the direction
of $h\in {\cal C}^1$ is given by
\begin{equation*}
\begin{split}
R'h(x\wedge y,&z\wedge u)=h(R(x,y)z,u)+g(R'h(x,y)z,u)\\
   &=h(R(x,y)z,u)+g((\nabla_y\nabla'h)(x,z),u)-g((\nabla_x\nabla'h)(y,z),u)\\
   &={1\over 2}\{\nabla^2_{yx}h(z,u)+\nabla^2_{yz}h(x,u)-\nabla^2_{yu}h(x,z)\\
  &\phantom{mmmm}- \nabla^2_{xy}h(z,u)-\nabla^2_{xz}h(y,u)+\nabla^2_{xu}h(y,z)
  \}+h(R(x,y)z,u)\\
&={1\over 2}\bigl\{ \nabla^2_{yz}h(x,u)+\nabla^2_{xu}h(y,z)-
\nabla^2_{xz}h(y,u)-\nabla^2_{yu}h(x,z)+\\
&\phantom{mmmmmmmm} h(R(x,y)z,u)-h(R(x,y)u,z)\bigr\}.
\end{split}
\end{equation*}
Where, in the last step, we have used the following identity
$$(\nabla^2_{xy}h-\nabla^2_{yx}h)(z,u)=h(R(x,y)u,z)+h(R(x,y)z,u).$$

Now using (\ref{DDh:form}), we get
$$R'h(x\wedge y,z\wedge u)= {1\over 2}\left\{-D\tilde Dh(x,y,z,u)
+h(R(x,y)z,u)-h(R(x,y)u,z)\right\}.$$
To get the derivative of $R$, as a symmetric curvature structure, i.e. in
${\cal C}^2$, it suffices to take the projection of the previous one, that
is
\begin{equation*}
\begin{split}
R'h(x\wedge y&,z\wedge u)={1\over 2}\biggl(R'h(x\wedge y,z\wedge u)+
R'h(z\wedge u,x\wedge y)\biggr)\\
&={1\over 4}\bigl\{-D\tilde Dh(x\wedge y,z\wedge u)-
D\tilde Dh(z\wedge u,x\wedge y)
+h(R(x,y)z,u)\\
&\phantom{mmmmmmm}-h(R(x,y)u,z)+
h(R(z,u)x,y)-h(R(u,z)x,y) \bigr\}.
\end{split}
\end{equation*}
This completes the proof of the  lemma. \ppp
\par\noindent\medskip
\begin{lemma} Let $(M,g)$ be a Riemannian manifold and
 $h\in {\cal C}^1$, then the differential of $h_{2k}$ at $g$, in the
direction of $h$,  is given by
$$h_{2k}'h={-1\over 2}<{c^{2k-1}\over (2k-1)!}R^k,h>-
  {k\over 4}(\delta\tilde\delta +\tilde\delta\delta )\biggl( *
  ( {g^{n-2k}\over (n-2k)!}R^{k-1}h)\biggr).$$
\end{lemma}
{\bf Proof.}
 Recall that
 $$h_{2k}=*({1\over (n-2k)!}g^{n-2k}R^k)= {1\over (n-2k)!}g^{n-2k}R^k(\mu_g,
 \mu_g),$$
 where $\mu_g$ is considered in the previous formula as an $n$-vector.
 then  using the previous lemma, we have
 \begin{equation*}
 \begin{split}
 h_{2k}'h=&{1\over (n-2k-1)!}g^{n-2k-1}hR^k(\mu_g,\mu_g)
 + {1\over (n-2k)!}g^{n-2k}kR^{k-1}R'h(\mu_g,\mu_g)\\
&{\phantom {mmmmmmmmmmm}}- 2 {1\over (n-2k)!}g^{n-2k}R^k({1\over 2}{\rm tr}_gh)
(\mu_g,\mu_g)\\
 =&*({1\over (n-2k-1)!}g^{n-2k-1}hR^k)
 + *({1\over (n-2k)!}g^{n-2k}kR^{k-1}R'h)\\
 &\phantom{mmmmmmmmmmmmmmmm}
- *( {1\over (n-2k)!}g^{n-2k}R^k){\rm tr}_gh\\
=&<T_{2k}-h_{2k}g,h>+ *\bigl\{ {g^{n-2k}kR^{k-1}\over 4(n-2k)!}
(-D\tilde Dh-\tilde D Dh+F_h(R)\bigr\}.\\
\end{split}
\end{equation*}
 Using first (\ref{dom:thet}) and then  (\ref{adjoint:dtilded}), we have
 \begin{equation*}
 \begin{split}
*\bigl( \frac{ g^{n-2k}kR^{k-1}}{ 4(n-2k)!}& (-D\tilde D h-\tilde D D h)\bigr)=
  *( D\tilde D +\tilde D D)(\frac{-kg^{n-2k}R^{k-1}h}{4(n-2k)!})\\
  &={-k\over 4}\{\delta\tilde\delta + \tilde\delta \delta \}\biggl( *
  ( {g^{n-2k}\over (n-2k)!}R^{k-1}h)\biggr).
\end{split}
\end{equation*}
On the other hand,  using simultaneously formulas
(\ref{fhomk}),(\ref{IPHodge}), (\ref{fh:autoadjoint}), proposition
\ref{extra} and formula (\ref{conjugg}), we get
 \begin{equation*}
 \begin{split}
*({1\over 4(n-2k)!})&g^{n-2k}kR^{k-1}(F_h(R))=
*({1\over 4(n-2k)!}g^{n-2k}F_h(R^k))\\
&={1\over 4}<F_h(R^k),*({1\over (n-2k)!}g^{n-2k})>
={1\over 4}<F_h(R^k),{1\over (2k)!}g^{2k}>\\
&={1\over 4}<R^k,F_h({1\over (2k)!}g^{2k})>
={1\over 4}<R^k,{2\over (2k-1)!}g^{2k-1}h>\\
&={1\over 2}<{1\over (2k-1)!}c^{2k-1}R^k,h>.
\end{split}
\end{equation*}
Recall that $T_{2k}=h_{2k}g-{1\over (2k-1)!}c^{2k-1}R^k$, then finally we have
$$h_{2k}'h=<-{1\over 2(2k-1)!}c^{2k-1}R^k,h>+
  {-k\over 4}(\delta\tilde\delta +\tilde\delta\delta )\biggl( *
  ( {g^{n-2k}\over (n-2k)!}R^{k-1}h)\biggr).$$
  The proof of the lemma is now complete. \ppp
  \par\smallskip\noindent
{\sc Remarks.} \begin{enumerate} \item In contrast with the volume form,  the derivative of $\mu_g$ seen
as an $n$-vector is $-{1\over 2}{\rm tr}_gh\mu_g$ and not just ${1\over 2}{\rm tr}_gh\mu_g$.
This fact was used
in the  previous proof.\\
\item The previous two lemmas are of local nature.
\end{enumerate}
\par\medskip\noindent
We are now ready to prove  the main
theorem. 
Note that for $k\in {\cal C}^1$,  both $\delta\tilde\delta k$ and 
$\tilde\delta\delta k$ are 
divergences
of some differential $1$-form. Therfore their integral is zero by Stokes' theorem.
 This applies particularly  to $k=*( {g^{n-2k}\over (n-2k)!}R^{k-1}h)\in {\cal C}^1$
of the previous formula, so that,
 \begin{equation*}
\begin{split}
H_{2k}'.h&=\int_M\biggl(h_{2k}'.h+{h_{2k}\over 2}{\rm
tr}_gh\biggr)\mu_g\\
&=-{1\over 2}<{c^{2k-1}\over (2k-1)!}R^k,h>+{h_{2k}\over 2}<g,h>\\
&={1\over 2}<h_{2k}g-{c^{2k-1}\over (2k-1)!}R^k,h>\\
&={1\over 2}<T_{2k},h>,
\end{split}
\end{equation*}
This completes the proof of the main theorem.
\ppp
\section{A Generalized Yamabe Problem}
As a consequence of the main theorem, we have  the following  result
\begin{proposition} For a compact Riemannian manifold $(M,g)$ with dimension
$n>2k$, the $(2k)$-Gauss-Bonnet curvature  $h_{2k}$ is constant if and only if $g$
is a critical point for the functional $H_{2k}$ when
restricted to the set ${\rm Conf}_0(g)$ of metrics pointwise conformal to $g$
and having the same total volume.\end{proposition}
{\bf Proof.} We proceed as in  the case of the scalar curvature \cite{Besse}.
 Note that $g$
is a critical point of $H_{2k}$ when restricted to ${\rm Conf}_0(g)$ if and
only if at $g$ we have
$$H_{2k}'.fg=0,$$
for all smooth functions $f$ such that $\int_Mf\mu_g=0$.\par\noindent
Next using  the main theorem, we have
\begin{equation*}
\begin{split}
H_{2k}'.fg=&<h_{2k}g-{c^{2k-1}R^k\over (2k-1)!},fg>\\
=&nfh_{2k}-f{c^{2k}R^k\over (2k-1)!}\\
=&(n-2k)fh_{2k}.
\end{split}
\end{equation*}
Consequently, for $\int_Mf\mu_g=0$, $H_{2k}'.fg=0$ if and only if the function $f$ is orthogonal
to $h_{2k}$.\par\noindent
Finally, consider the function
$$f=h_{2k}-{1\over {\rm vol}(g)}\int_Mh_{2k}\mu_g,$$
it is orthogonal to $h_{2k}$ and to the constants, then it is the zero
function. This completes the proof.\ppp
\par\medskip\noindent
It is then natural  to ask whether for each $k$ we have:\par\medskip\noindent
{\bf Question:} ({\bf Generalized Yamabe problem.})\par\smallskip\noindent
{\sl In  each conformal class of
a fixed Riemannian metric on a smooth compact manifold with dimension $n>2k$
 there exists a metric with $h_{2k}$ constant.}
 \par\medskip\noindent
 A closely related and at the same time  parallel problem to the previous one is the 
 $\sigma_k$-Yamabe problem. It is at present the
subject of intensive researches. We can state it as follows: \\
{\sl Let   $\sigma_k$ denote the symmetric function of order $k$ in the eigenvalues of the 
 Schouten tensor. For each $k$, there exists 
in  each conformal class of
a fixed Riemannian metric on a smooth compact manifold 
  a metric with $\sigma_k$ constant.}\\
\par\smallskip
These two problems coincide in the class of a conformally flat metric and when 
 $k$ is even. In fact, in this case the curvatures
  $h_{2k}$ and  $\sigma_{2k}$ differ only by a constant factor, \cite{Labbi}.
Remark also, that both problems generalize the classical Yamabe problem obtained for $k=1$ and
$k=2$ respectively.  \\
The $\sigma_k$-Yamabe problem was recently proved for $k>n/2$  by
Gursky and Viaclovsky
\cite{Gursky} with the assumption that the original metric is  ``admissible''.
Then  Sheng, Trudinger and
Wang \cite{Trudinger} completed the proof of the remaining cases where $2\leq k\leq n/2$ 
after assuming that 
the relevant equation  is variational. 
\\
Note also that the $\sigma_k$-Yamabe  problem was  solved  in the conformally flat case by
  Li-Li \cite{Lili}, and   Guan-Wang \cite{Guan}.\\

{\sc Remark.} 
The curvatures  $h_{2k}$ are in general different from the symmetric functions 
in the eigenvalues of the Riemann curvature operator $R$. In fact, they  are not even in general
symmetric polynomials in the eigenvalues of $R$ as one can check it easily. In fact, for
a $4$-dimensional Riemannian manifold,  suppose
the eigenvectors of $R$ are all of rank $1$ (decomposed)  and let $\lambda_1,...,\lambda_6$
denote the eigenvalues of $R$. Let us re-arrange  so that $\lambda_1, \lambda_2$
(resp. $\lambda_3, \lambda_4$ and $\lambda_5, \lambda_6$) are two eigenvalues corresponding
to supplementary $2$-planes (eigenvectors), then $h_4$, up to a constant, equals 
$\lambda_1\lambda_2+\lambda_3\lambda_4+\lambda_5\lambda_6$. This is clearly not a symmetric 
polynomial in $\lambda_1,...,\lambda_6$.

\section{Generalized Einstein Manifolds}
The usual Einstein metrics are the critical metrics for the total scalar curvature functional
when restricted to those metrics with  unit volume. Similarly, considering the critical
metrics of the functional associated to the 
 Gauss-Bonnet curvatures we obtain generalized Einstein metrics. Precisely:\\ 

 \begin{definition}
For $2\leq 2k\leq n$, we say that $(M,g)$ is $(2k)$-Einstein if its $(2k)$-Einstein-Lovelock tensor
is proportional to the metric, that is
$$T_{2k}=\lambda g.$$
 \end{definition}
 Note that since the tensors $T_{2k}$ are divergence free, the function $\lambda$
 is then a constant.\par\noindent
 Also, the $2$-Einstein manifolds are the usual Einstein manifolds,
 and when $n=2k$, we have $T_{2k}=0$ for any metric on
 $M.$\par\medskip\noindent
 The class of ($2k$)-Einstein Riemannian manifolds contains the
 manifolds with constant sectional curvature, and all isotropy
 irreducible homogeneous manifolds with their canonical Riemannian
 metrics.  
\par\medskip\noindent
{\sc Examples.}
\begin{itemize}
\item  Let M be an arbitrary Riemannian manifold with dimension
$n\geq 3$, the  Riemannian product $M\times {\bf R}^q$ is with $T_{2k}=0$
for $2k\geq n$ but $T_2$ is not proportional to the metric. This example shows
that in some sens the generalized Einstein condition is weaker
then the usual one. But this is not always true, as shown by the
next example:
\item Let $M$ be a 4-dimensional Ricci-flat but not flat manifold
(for example a K3 surface endowed with the Calabi-Yau metric), and
$T^q$ be a flat torus. Then the Riemannian product $M\times T^q$
is Einstein in the usual sens ($T_2=0$), but it is not 4-Einstein.
\item Let $(M,g)$ be a conformally flat manifold, then it can be
shown without difficulties that if $(M,g)$ is ($2k$)-Einstein then the Ricci tensor has
at each point at  most k distinct eigenvalues. Similar results
hold for hypersurfaces of the Euclidean space.
\end{itemize}
\par\medskip\noindent

 With respect to the orthogonal decomposition into irreducible components of the ring of curvature structures \cite{Kulk, Labbi},
  the Gauss-Kronecker tensor $R^q$ splits to
  \begin{equation}\label{decomp}
R^q=\omega_{2q}+g\omega_{2q-1}+...+g^{2q-1}\omega_{1}+g^{2q}\omega_0.
\end{equation}

Evidently,  the tensor $T_{2k}$ is proportional to the metric if and only 
if the
generalized Ricci tensor $c^{2k-1}R^k$ is also proportional to the metric. Therefore
it results immediately from lemma 5.7 of \cite{Labbi} the following proposition.
\begin{proposition}
A Riemannian metric is ($2q$)-Einstein if and only if the component
$\omega_1$ of $R^q$ with respect to the irreducible splitting (\ref{decomp}) vanishes.
\end{proposition}
The previous result  generalizes a similar well known  result about the usual Einstein
 metrics.\par\medskip\noindent

{\bf Acknowledgements.} I would like to thank J. Lafontaine  and R. Souam 
 for useful discussions.

\noindent
Mohammed-Larbi Labbi\\
 Department of Mathematics,\\
  College of Science, University of Bahrain,\\
  P. O. Box  32038,  Bahrain.\\
  E-mail: labbi@sci.uob.bh

\end{document}